\documentstyle{amsppt}
\def\ve{\varepsilon}
\topmatter
\author{S.~V.~Konyagin}
\endauthor
\title{A sum-product estimate in fields of prime order}
\endtitle
\abstract
Let $q$ be a prime, $A$ be a subset of a finite field $F:=\Bbb Z/q\Bbb Z$,
$|A|<\sqrt{|F|}$. We prove the estimate 
$\max(|A+A|,|A\cdot A|)\ge c|A|^{1+\ve}$ for some $\ve>0$ and $c>0$.
This extends the result of \cite{BKT}.
\endabstract
\endtopmatter
\document
Key words: subsets of finite fields, groups.

MSC: 11B75, 11T30.

\centerline{\S1. Introduction}
\bigskip
Let $q$ be a prime, $F=\Bbb Z/q\Bbb Z$, $F^*=F\setminus\{0\}$, and let $A$ be 
a nonempty subset of $F$. We consider the sum set
$$A+A:=\{a+b:a,b\in A\}$$
and the product set
$$A\cdot A:=\{ab:a,b\in A\}.$$
Let $|A|$ denote the cardinality of $A$. We have the obvious bounds
$$|A+A|,|A\cdot A|\ge|A|.$$
The bounds are clearly sharp if $A=F$ or if $|A|=1$. We can expect 
some improvement in other cases. However, good lower estimates 
for $\max(|A+A|,|A\cdot A|)$ were not known for a long time.
Recently a breakthrough was made by J.~Bourgain, N.~Katz, and
T.~Tao\cite{BKT} who proved the following result.
\proclaim{Theorem A} Let $A$ be a subset of $F$ such that
$$|F|^\delta<|A|<|F|^{1-\delta}$$
for some $\delta>0$. Then one has a bound of the form
$$\max(|A+A|,|A\cdot A|)\ge c(\delta)|A|^{1+\ve}$$
for some $\ve=\ve(\delta)>0$ and $c(\delta)>0$.
\endproclaim
Also, in \cite{BKP} a reader can find the history of the problem, 
generalizations and applications of Theorem A.

However, Theorem A does not estimate $\max(|A+A|,|A\cdot A|)$ if
$|A|$ is small comparatively to $|F|$. The aim of this paper is to 
give such an estimate.
\proclaim{Theorem 1} Let $A$ be a subset of $F$ such that
$$|A|<|F|^{1/2}.$$
Then one has a bound of the form
$$\max(|A+A|,|A\cdot A|)\ge c|A|^{1+\ve}$$
for some $\ve>0$ and $c>0$.
\endproclaim
Clearly, Theorem A and Theorem 1 immediately imply uniform estimates for
$|A|<|F|^{1-\delta}$.
\proclaim{Corollary 1} Let $A$ be a subset of $F$ such that
$$|A|<|F|^{1-\delta}$$
for some $\delta>0$. Then one has a bound of the form
$$\max(|A+A|,|A\cdot A|)\ge c(\delta)|A|^{1+\ve}$$
for some $\ve=\ve(\delta)>0$ and $c(\delta)>0$.
\endproclaim
To prove Theorem A, the authors associated with a set $A\subset F$
the following set
$$I(A):=\{a_1(a_2-a_3)+a_4(a_5-a_6):a_1,\dots,a_6\in A\}.$$
They found lower bounds for $|I(A)|$ and applied those bounds  
for estimation of $\max(|A+A|,|A\cdot A|)$. Using the main idea of
\cite{BKT} we give new lower estimates for $|I(A)|$.

We denote
$$A-A:=\{a-b:a,b\in A\}.$$
Throughout the paper $c$ and $C$ will denote absolute positive constants.
If $f$ and $g$ are functions, we will write $f\ll g$ or $g\gg f$ if 
$|A|\le CB$ for some constant $C$ (uniformly with respect to variables in
$f$ and $g$).
\proclaim{Theorem 2} Let $A$ be a subset of $F$ such that
$$|A|<|F|^{1/2}.$$
Then one has a bound of the form
$$|A-A|\times|I(A)|\ge c|A|^{5/2}.$$
\endproclaim
Observing that $|I(A)|\ge|A-A|$ we deduce from Theorem 2 an estimate for
$|I(A)|$.
\proclaim{Corollary 2} Let $A$ be a subset of $F$ such that
$$|A|<|F|^{1/2}.$$
Then one has a bound of the form
$$|I(A)|\ge c|A|^{5/4}.$$
\endproclaim
Also, one can get a good lower bound for $|I(A)|$ if $|A|>|F|^{1/2}$.
\proclaim{Theorem 3} Let $A$ be a subset of $F$ such that
$$|A|>|F|^{1/2}.$$
Then one has a bound of the form
$$|I(A)|\ge|F|/2.$$
\endproclaim
\bigskip
\centerline{\S2. A version of a lemma of J.~Bourgain, N.~Katz, and
T.~Tao} 
\centerline{and the proof of Theorem 3}
\bigskip
Let $A$ be a subset of $F$ and $\xi\in F$. Denote
$$S_\xi(A):=\{a+b\xi:a,b\in A\}.$$
We will use the following analog of Lemma 4.2 from \cite{BKT}.
\proclaim{Lemma 1} Let $\xi\in F^*$ and 
$$|S_\xi(A)|<|A|^2.$$
Then 
$$|I(A)|\ge|S_\xi(A)|.$$
\endproclaim
\demo{Proof} By the supposition on $|S_\xi(A)|$, the surjection
$$A^2\to F:\,(a,b)\to a+b\xi$$
cannot be one-to-one. Thus there are $(a_1,b_1)\neq(a_2,b_2)$
with
$$(a_1-a_2)+(b_1-b_2)\xi=0.\tag1$$ 
We observe that $b_1\neq b_2$. Denote
$$S:=(b_1-b_2)S_\xi(A):=\{(b_1-b_2)s:s\in S_\xi(A)\}.$$
We have 
$$|S|=|S_\xi(A)|,$$
and every element $s\in S$ can be represented in a form
$$s=(b_1-b_2)a+(b_1-b_2)b\xi,\quad a,b\in A.$$
Substituting $(b_1-b_2)\xi$ from (1) we get
$$s=(b_1-b_2)a+(a_2-a_1)b.$$
Therefore, $S\subset I(A)$, and Lemma 1 is proved.
\enddemo

To prove Theorem 3 with a slightly weaker form we can use Lemma 2.1 from 
\cite{BKP}. However, we shall need a following generalization of that lemma.
\proclaim{Lemma 2} Let $A\subset F$, $G\subset F^*$. Then there exists 
$\xi\in G$ such that
$$|S_\xi(A)|\ge|A|^2|G|/(|A|^2+|G|).$$
\endproclaim
\demo{Proof} Denote for $\xi\in G$ and $s\in F$
$$f_\xi(s):=\{(a,b):a,b\in A,a+b\xi=s\}.$$
We have
$$\gather
\sum_{s\in F}f_\xi(s)^2=|\{(a_1,b_1,a_2,b_2):a_1,b_1,a_2,b_2\in A,
a_1+b_1\xi=a_2+b_2\xi\}|\\
=|A|^2+|\{(a_1,b_1,a_2,b_2):a_1,b_1,a_2,b_2\in A,a_1\neq a_2,
a_1+b_1\xi=a_2+b_2\xi\}|.
\endgather$$
Taking the sum over $\xi\in G$ and observing that for any 
$a_1,b_1,a_2,b_2\in A$ with $a_1\neq a_2$ there is at most one
$\xi$ such that $a_1+b_1\xi=a_2+b_2\xi$ we obtain
$$\sum_{g\in G}\sum_{s\in F}f_\xi(s)^2\le|A|^2|G|+|A|^4.$$
Therefore, we can fix $\xi\in G$ so that
$$\sum_{s\in F}f_\xi(s)^2\le|A|^2+|A|^4/|G|.\tag2$$
(Clearly we can assume that $G\neq\emptyset$.) 
By Cauchy---Schwartz inequality,
$$\left(\sum_{s\in F}f_\xi(s)\right)^2\le|S_\xi(A)|\sum_{s\in F}f_\xi(s)^2.$$
Moreover,
$$\sum_{s\in F}f_\xi(s)=|A|^2.$$
Therefore, by (2),
$$|S_\xi(A)|\ge|A|^4/(|A|^2+|A|^4/|G|)=|A|^2|G|/(|A|^2+|G|).$$
\enddemo

\demo{The proof of Theorem 3} Take $G=F^*$. 
Taking into account the supposition on $|A|$ we deduce from
Lemma 2 that for some $\xi$
$$|S_\xi(A)|\ge|A|^2|G|/(|A|^2+|G|)>|A|^2|G|/(2|A|^2)=(|F|-1)/2$$
Therefore, $|S_\xi(A)|\ge|F|/2$. Also, we have
$|S_\xi(A)|\le p<|A|^2$. Thus, Theorem 3 follows from Lemma 1.
\enddemo
\bigskip
\centerline{\S3. Some preparations}
\bigskip
Let $A\subset F^*$ and 
$$H:=\{s\in F:|\{(a,b):a,b\in A, s=a/b\}|\ge|A|^2/(5|A\cdot A|)\}.$$
Denote by $G$ the multiplicative subgroup of $F^*$ generating by $H$.
\proclaim{Lemma 3} There is a coset $G_1$ of $G$
such that 
$$|A\cap G_1|\ge|A|/3.\tag3$$
\endproclaim
\demo{Proof} Assume the contrary. Let $A_1, A_2,\dots$ be the nonempty 
intersections of $A$ with cosets of $G$. Take a minimal $k$ so that
$$\left|\bigcup_{i=1}^k A_i\right|\ge|A|/3$$
and denote
$$A'=\bigcup_{i=1}^k A_i,\quad A''=A\setminus A'.$$
We have
$$|A'|\ge|A|/3.$$
On the other hand,
$$|A'|\le\left|\bigcup_{i=1}^{k-1} A_i\right|+|A_k|<2|A|/3.$$
Hence,
$$|A|/3\le|A'|\le2|A|/3$$
and
$$|A'|\times|A''|=|A'|(|A|-|A'|)\ge2|A|^2/9.\tag4$$
Denote for $s\in F^*$
$$f(s):=\{(a,b):a\in A',\,b\in A'',a/b=s\}.$$
Note that if $a\in A',\,b\in A'',$ then $a/b\not\in H$. Therefore,
for any $s$ we have the inequality $f(s)<|A|^2/(5|A\cdot A|)$.
Hence, 
$$
\sum_{s\in F^*}f(s)^2\le\frac{|A|^2}{5|A\cdot A|}\sum_{s\in F^*}f(s)
=\frac{|A|^2|A'|\times|A''|}{5|A\cdot A|}.\tag5
$$
Denote for $s\in F^*$
$$g(s):=\{(a,b):a\in A',\,b\in A'',ab=s\}.$$
By Cauchy---Schwartz inequality,
$$\left(\sum_{s\in F}g(s)\right)^2\le|A\cdot A|\sum_{s\in F}g(s)^2.$$
Therefore,
$$\sum_{s\in F^*}g(s)^2\ge\left(\sum_{s\in F}g(s)\right)^2/|A\cdot A|
=\frac{(|A'|\times|A''|)^2}{|A\cdot A|}.\tag6$$
Now observe that both the sums $\sum_{s\in F^*}f(s)^2$ and 
$\sum_{s\in F^*}g(s)^2$ are equal to the number of solutions of the equation
$a_1'a_1''=a_2'a_2''$, $a_1',a_2'\in A'$, $a_1'',a_2''\in A''$. Thus,
comparing (5) and (6) we get
$$|A'|\times|A''|\le|A|^2/5.$$
But the last inequality does not agree with (4), and the proof is complete.
\enddemo

We will use the function $S_\xi(A)$ defined in the beginning of \S2.
\proclaim{Lemma 4} Let $A\subset F^*$ and $|A|>1$. Then there exists 
$\xi\in G$ such that
$$\min\left(|A|^3/(5|A\cdot A|),|A|^2|G|/(|A|^2+|G|)\right)\le S_\xi(A)
<|A|^2.\tag7$$
\endproclaim
\demo{Proof} We consider two cases.

1. Case 1: there exists $g\in G$ such that $S_g(A)=|A|^2$. We claim that  
$$\exists\xi\in G\quad |A|^3/(5|A\cdot A|)\le S_\xi(A)<|A|^2.\tag8$$
Assume that (8) does not hold. Take an arbitrary $g\in G$ satisfying 
$S_g(A)=|A|^2$ (this means that the elements $a+bg$, $a,b\in A$
are pairwise distinct) and an arbitrary $h\in H$. Denote 
$$A_h=\{b\in A: bh\in A\}.$$
We have $|A_h|\ge|A|^2/(5|A\cdot A|)$ because $h\in H$. By our supposition on
$g$, all the sums $a+b(gh)=a+(bh)g$, $a\in A$, $b\in A_h$, are distinct.
Therefore, $S_{gh}(A)\ge|A|^3/(5|A\cdot A|)$. Our supposition that (8) 
does not hold implies that $S_{gh}(A)=|A|^2$.

So, we see that if an elements $g\in G$ satisfies the condition 
$S_g(A)=|A|^2$ then for any $h\in H$ the elements $gh$ also satisfies this 
condition. Since $H$ generates $G$ we deduce that the condition
$S_g(A)=|A|^2$ holds for all elements $g\in G$. But this is impossible
because $S_1(A)\le|A|(1+|A|)/2<|A|^2$, and (8) is proved. 

2. Case 2: for all $g\in G$ we have $S_g(A)<|A|^2$. Then the existence of 
a required $\xi\in G$ immediately follows from Lemma 2, and the proof of 
Lemma 4 is complete.  
\enddemo

\proclaim{Lemma 5} Let $G$ be a subgroup of $F^*$, $B\subset G$, 
$|B|<\sqrt{|F|}$. Then
$$|B-B|\gg|A|^{5/2}/|G|.$$
\endproclaim
\demo{Proof} We use arguments from \cite{HBK}. Consider the cosets 
$G_1, G_2,\dots$ of $G$ in $F^*$. For any coset $G_t$ and $s\in G_t$ denote
$$N_t:=|\{(g_1,g_2):\,g_1,g_2\in G, g_1-g_2=s\}|.$$
(It is clear that the definition is correct, namely, it does not depend on 
the choice of $s\in G_t$.) We order $G_1, G_2,\dots$ in such a way that
$$N_1\ge N_2\dots.\tag9$$
Also, we consider $N_t=0$ if $t$ exceeds the number of the cosets of $G$
in $F^*$. Lemma 5 from \cite{HBK} claims that if 
$$|G|^4T<|F|^3,\tag10$$
then
$$\sum_{t=1}^T N_t\ll(|G|T)^{2/3}.\tag11$$
Therefore, if (10) holds, then, by (9) and (11),
$$N_T\ll |G|^{2/3}T^{-1/3},\tag12$$
and, moreover,
$$\sum_{t=1}^T N_t^2\ll \sum_{t=1}^T \left(|G|^{2/3}t^{-1/3}\right)^2
\ll |G|^{4/3}T^{1/3}.\tag13$$

Now denote for $t=1,2,\dots$ 
$$L_t:=|\{(b_1,b_2):\,b_1,b_2\in B, b_1-b_2\in G_t\}|,$$
$$M_t:=|\{(b_1,b_2,b_3,b_4):\,b_1,b_2,b_3,b_4\in B, 
b_1-b_2=b_3-b_4\in G_t\}|.$$
We consider $L_t=M_t=0$ if $t$ exceeds the number of the cosets of $G$.
For every element $b_1\in B$ there are at most $N_t$ elements $b_2\in B$
such that $b_1-b_2\in G_t$. Therefore,
$$L_t\le N_t|B|.\tag14$$
Further, for every element $s\in G_t$ there are at most $N_t$ elements 
$b_1\in B$ such that $b_1-s\in B$. Therefore,
$$M_t\le N_tL_t\le N_t^2|B|.\tag15$$
Also,
$$\sum_t L_t=|\{(b_1,b_2):\,b_1,b_2\in B, b_1\neq b_2\}|=|B|(|B|-1).\tag16$$

The statement of Lemma 5 is trivial if $|G|\ge|B|^{3/2}$. Thus, we assume that
$$|G|<|B|^{3/2}.\tag17$$
Let us take
$$T:=[|B|^{3/2}/|G|].$$
By the supposition on $|B|$ and (17), we have
$$|G|^4T\le |G|^3|B|^{3/2}<|B|^6<|F|^3,$$
and (10) holds. By (15) and (13),
$$\sum_{t=1}^T M_t\le |B|\sum_{t=1}^T N_t^2\ll |B||G|^{4/3}T^{1/3}.$$
By (15), (12), and (16),
$$\sum_{t>T} M_t\le\sum_{t>T} N_tL_t\le N_T\sum_{t}L_t\ll|G|^{2/3}T^{-1/3}
|B|^2.$$
Thus, taking into account the definition of $T$, we get
$$\sum_t M_t\ll|B||G|^{4/3}T^{1/3}+|B|^2|G|^{2/3}T^{-1/3}
\ll |B|^{3/2}|G|.\tag18$$

Denote for $s\in F^*$
$$f(s):=|\{(b_1,b_2):b_1,b_2\in B,b_1-b_2=s\}|.$$
By Cauchy---Schwartz inequality,
$$\left(\sum_{s\in F^*}f(s)\right)^2\le\sum_{s\in F^*}f(s)^2
|\{b_1-b_2: b_1,b_2\in G,b_1\neq b_2\}|.$$
Also,
$$\sum_{s\in F^*}f(s)=\sum_t L_t,\quad \sum_{s\in F^*}f(s)^2=\sum_t M_t,$$
and, by (16) and (18), we find
$$|B-B|\gg 1+(|B|(|B|-1))^2/(|B|^{3/2}|G|)\gg |B|^{5/2}/|G|.$$
Lemma 5 is proved.
\enddemo
\bigskip
\centerline{\S4. The proofs of Theorems 1 and 2}
\bigskip
\demo{Proof of Theorem 2} In the case $|A\cdot A|>|A|^{3/2}$ the assertion
is trivial because $|A-A|\ge|A|$, $|I(A)|\ge|A\cdot A|$. Thus, we will 
consider that $|A\cdot A|\le|A|^{3/2}$ and, therefore, 
$$|A|^3/(|A\cdot A|)\ge|A|^{3/2}.\tag19$$
We take the group $G$ defined in the beginning of \S3. If $|G|>|A|^{3/2}$,
then
$$|A|^2|G|/(|A|^2+|G|)\ge|A|^2|A|^{3/2}/(|A|^2+|A|^{3/2})\ge|A|^{3/2}/2.$$
Hence, by Lemma 4 and (19), there exists $\xi\in G$ such that
$$|A|^{3/2}\ll S_\xi(A)<|A|^2.$$
Lemma 1 claims that $|I(A)|\gg|A|^{3/2}$. Thus,
$$|A-A|\times|I(A)|\gg|A|\times|A|^{3/2},$$
and we have the required estimate. 

Now it suffices to consider the case if both the conditions (19) and
$$|G|\le|A|^{3/2}\tag20$$
are satisfied. We see from (20) that
$$|A|^2|G|/(|A|^2+|G|)\ge|A|^2|G|/(2|A|^2)=|G|/2.$$
Hence, by Lemma 4 and (19), there exists $\xi\in G$ such that
$$|G|\ll S_\xi(A)<|A|^2,$$
and, by Lemma 1,
$$|I(A)|\gg|G|.\tag21$$

We take a coset $G_1$ of $G$ in accordance with Lemma 3. Fix an arbitrary
$g_1\in G_1$. Let
$$B:=\{b\in G: bg_1\in A\}.$$
By (3), we have 
$$|B|=|A\cap G_1|\ge|A|/3.$$
We get from Lemma 5 that
$$|A-A|\ge|A\cap G_1-A\cap G_1|=|B-B|\gg(|A|/3)^{5/2}/|G|\gg|A|^{5/2}/|G|.
\tag22$$
Now Theorem 2 follows from (21) and (22).
\enddemo

Note that $|A-A|\le|I(A)|$. Thus, Theorem 2 implies the inequality 
$$|I(A)|\gg|A|^{5/4}\tag23$$
provided that $|A|<|F|^{1/2}$. Theorem 1 is a corollary of (23) and Lemma 2.4 
from \cite{BKT}.

The research was fulfilled during the author's visit to the University Aroma 
ARE and was supported by the Mathematical Department of the University
and the Institute Nazionale di Alta Matematica.

\bigskip
\centerline{References}
\bigskip
[BKT] J.~Bourgain, N.~Katz, and T.~Tao, A sum-product estimate in finite fields
and their applications. ArXiv:math.CO/0301343 v1, January 29, 2003.

[HBK] D.~R.~Heath-Brown, S.~V.~Konyagin, New bounds for Gauss sums derived 
from $k$th powers, and for Heilbronn's exponential sums. Quart.~J.~Math.
{\bf 51} (2000), 221--235.

Department of Mechanics and Mathematics, Moscow State University,\newline 
Moscow, 119992, Russia.

E-mail address: konyagin\@ok.ru

\enddocument